\documentclass[a4paper,11pt,reqno]{elsarticle}
\usepackage{latexsym, amsmath, amsfonts, amssymb, amsthm, amscd,epsfig}
\usepackage{url}

\usepackage[font=small, labelfont={sc}]{caption}
\usepackage{subfigure}

\usepackage[a4paper,scale={0.73,0.75},marginratio={1:1},footskip=10mm,headsep=10mm]{geometry}

\usepackage{color}

\usepackage{hyperref}

\numberwithin{equation}{section}

\def \beq {\begin{eqnarray}}
\def \eeq {\end{eqnarray}}
\def \beqn {\begin{eqnarray*}}
\def \eeqn {\end{eqnarray*}}

\newtheorem{theorem}{Theorem}
\newtheorem{itlemma}[theorem]{Lemma}
\newtheorem{itproposition}[theorem]{Proposition}
\newtheorem{itcorollary}[theorem]{Corollary}
\newtheorem{itremark}[theorem]{Remark}
\newtheorem{itdefinition}[theorem]{Definition}
\newtheorem{itexample}[theorem]{Example}
\newtheorem{itclaim}[theorem]{Claim}
\newtheorem{itfact}[theorem]{Fact}

\newenvironment{fact}{\begin{itfact}\rm}{\end{itfact}}
\newenvironment{claim}{\begin{itclaim}\rm}{\end{itclaim}}
\newenvironment{lemma}{\begin{itlemma}}{\end{itlemma}}
\newenvironment{remark}{\begin{itremark}\rm}{\end{itremark}}
\newenvironment{corollary}{\begin{itcorollary}}{\end{itcorollary}}
\newenvironment{proposition}{\begin{itproposition}}{\end{itproposition}}
\newenvironment{definition}{\begin{itdefinition}\rm}{\end{itdefinition}}
\newenvironment{example}{\begin{itexample}\rm}{\end{itexample}}

\newcommand{\be}[1]{\begin{equation}\label{#1}}
\newcommand{\ee}{\end{equation}}
\newcommand{\bl}[1]{\begin{lemma}\label{#1}}
\newcommand{\br}[1]{\begin{remark}\label{#1}}
\newcommand{\brs}[1]{\begin{remarks}\label{#1}}
\newcommand{\bt}[1]{\begin{theorem}\label{#1}}
\newcommand{\bd}[1]{\begin{definition}\label{#1}}
\newcommand{\bp}[1]{\begin{proposition}\label{#1}}
\newcommand{\bc}[1]{\begin{corollary}\label{#1}}
\newcommand{\bfact}[1]{\begin{fact}\label{#1}.}
\newcommand{\bex}[1]{\begin{example}\label{#1}.}
\newcommand{\ec}{\end{corollary}}
\newcommand{\efact}{\end{fact}}
\newcommand{\eex}{\end{example}}
\newcommand{\el}{\end{lemma}}
\newcommand{\er}{\end{remark}}
\newcommand{\ers}{\end{remarks}}
\newcommand{\et}{\end{theorem}}
\newcommand{\ed}{\end{definition}}
\newcommand{\ep}{\end{proposition}}
\newcommand{\epr}{\end{proof}}
\newcommand{\bpr}{\begin{proof}}
\newcommand{\bcl}[1]{\begin{claim}\label{#1}}
\newcommand{\ecl}{\end{claim}}

\newcommand{\ecs}{\end{corollary}}
\newcommand{\eers}{\end{exercise}}
\newcommand{\eexs}{\end{example}}
\newcommand{\eems}{\end{example}}
\newcommand{\els}{\end{lemma}}
\newcommand{\eles}{\end{lemmaex}}
\newcommand{\ets}{\end{theorem}}
\newcommand{\eds}{\end{definition}}
\newcommand{\eps}{\end{proposition}}

\newcommand{\bi}{\begin{itemize}}
\newcommand{\ei}{\end{itemize}}
\newcommand{\ben}{\begin{enumerate}}
\newcommand{\een}{\end{enumerate}}

\def\vbar{\mathchoice{\vrule height6.3ptdepth-.5ptwidth.8pt\kern-.8pt}
   {\vrule height6.3ptdepth-.5ptwidth.8pt\kern-.8pt}
   {\vrule height4.1ptdepth-.35ptwidth.6pt\kern-.6pt}
   {\vrule height3.1ptdepth-.25ptwidth.5pt\kern-.5pt}}
\def\fudge{\mathchoice{}{}{\mkern.5mu}{\mkern.8mu}}
\def\bbc#1#2{{\rm \mkern#2mu\vbar\mkern-#2mu#1}}
\def\bbb#1{{\rm I\mkern-3.5mu #1}}
\def\bba#1#2{{\rm #1\mkern-#2mu\fudge #1}}
\def\bb#1{{\count4=`#1 \advance\count4by-64 \ifcase\count4\or\bba A{11.5}\or
   \bbb B\or\bbc C{5}\or\bbb D\or\bbb E\or\bbb F \or\bbc G{5}\or\bbb H\or
   \bbb I\or\bbc J{3}\or\bbb K\or\bbb L \or\bbb M\or\bbb N\or\bbc O{5} \or
   \bbb P\or\bbc Q{5}\or\bbb R\or\bbc S{4.2}\or\bba T{10.5}\or\bbc U{5}\or
   \bba V{12}\or\bba W{16.5}\or\bba X{11}\or\bba Y{11.7}\or\bba Z{7.5}\fi}}

\def \R {{\mathbb R}}

\def \E {{\mathbb E}}

\def \ra {\rightarrow }
\def \x {\underline{x}}

\def \s {y}

\newcommand{\ba}[1]{\addtocounter{for}{1} \begin{eqnarray}\label{#1}}
\newcommand{\ea}{\end{eqnarray}}
\def\Box{{\hfill\hbox{\enspace${\square}$}} \smallskip}
\def\sqr#1#2{{\vcenter{\vbox{\hrule height .#2pt
                             \hbox{\vrule width .#2pt height#1pt \kern#1pt
                                   \vrule width .#2pt}
                             \hrule height .#2pt}}}}
\def\square{\mathchoice\sqr54\sqr54\sqr{4.1}3\sqr{3.5}3}
\def\pmb#1{\setbox0=\hbox{#1}%
   \kern-.025em\copy0\kern-\wd0
   \kern.05em\copy0\kern-\wd0
   \kern-.025em\raise.0433em\box0 }
\def\sqr#1#2{{\vcenter{\vbox{\hrule height.#2pt
     \hbox{\vrule width.#2pt height#1pt \kern#1pt
   \vrule width.#2pt}\hrule height.#2pt}}}}

\def\s{\sigma}
\def\d{\delta}
\def\l{\lambda}

\def\a{\alpha}
\def\b{\beta}

\newenvironment{myenumerate}{%
\begin{list}{\labelenumi}
	{%
	\setlength{\itemsep}{0.4em}%
	\setlength{\topsep}{0.5em}%
	\setlength\leftmargin{2.6em}%
	\setlength\labelwidth{2.15em}%
	\setlength{\labelsep}{0.45em}%
	\usecounter{enumi}%
	}%
	}%
{\end{list}
}

\renewenvironment{enumerate}{
\renewcommand{\theenumi}{\arabic{enumi}}%
\renewcommand{\labelenumi}{{\rm(\theenumi)}}%
\begin{myenumerate}}%
{\end{myenumerate}}

{\end{myenumerate}}

{\end{myenumerate}}

\newenvironment{myitemize}{%
\begin{list}{$\bullet$}%
 	{%
	\setlength{\itemsep}{0.4em}%
	\setlength{\topsep}{0.5em}%
	\setlength\leftmargin{2.6em}%
	\setlength\labelwidth{2.15em}%
	\setlength{\labelsep}{0.45em}%
	}%
	}%
{\end{list}}

\renewenvironment{itemize}{
\begin{myitemize}}%
{\end{myitemize}}

\begin{document}

\title{Collective periodicity in mean-field models of cooperative behavior}

\author[fc]{F. Collet}
\ead{francesca.collet@unibo.it}

\author[pdp]{P. Dai Pra
}
\ead{daipra@math.unipd.it}

\author[mf]{M. Formentin}
\ead{formentin@utia.cas.cz}

\address[fc]{Dipartimento di Matematica, Alma Mater Studiorum Universit\`a di Bologna, Piazza di Porta San Donato~5, I-40126 Bologna, Italy}

\address[pdp]{Dipartimento di Matematica, Universit\`a degli Studi di Padova, Via Trieste 63, I-35121 Padova, Italy}

\address[mf]{UTIA, Czech Academy of Sciences, Pod Vod\'{a}renskou v\v{e}\v{z}\'{i} 4, CZ-18208 Prague, Czech Republic}

%

\date{\today}

\begin{abstract} 
We propose a way to break symmetry in stochastic dynamics by introducing a dissipation term. We show in a specific mean-field model, that if the reversible model undergoes a phase transition of ferromagnetic type, then its dissipative counterpart exhibits periodic orbits in the thermodynamic limit.

\end{abstract}

\maketitle

\section{Introduction}

Many real systems comprised by many interacting components, for example neuronal networks, may exhibit collective periodic behavior even though single components have no ``natural'' tendency of behaving periodically. Various {\em stylized} models have been proposed to capture the essence of this phenomenon; we refer the reader to \cite{lens} and \cite{gp}, where many related references are given, and possible origins of this periodic behavior are discussed in clear and rigorous terms. Most of the proposed models, despite of their rather simple structure, turn out to be hard to analyze in rigorous terms; their study ends up to the search of stable attractors of nonlinear, infinite dimensional dynamical systems. 

The main purpose of this paper is to illustrate and study a class of models that are obtained as perturbations of ``classical'' reversible models by adding a dissipation term, which breaks reversibility. Models of this sort have already appeared in different context (see e.g. \cite{Schweitzer}  for applications to multicellular dynamics, and \cite{gsss} for applications to finance). The simplest interacting system in this class is the dissipative Curie-Weiss model proposed and analyzed in \cite{dfr}. Its infinite volume limit can be reduced to a two dimensional nonlinear dynamics, whose large-time behavior can be determined explicitly. It is shown, in particular, that by tuning the parameters of the model a globally stable periodic orbit emerges through a {\em Hopf bifurcation}. Note that this phenomenon is qualitatively similar to the emergence of periodic orbit in a system of interacting FitzHugh-Nagumo models presented in \cite{lens}.

In this paper, after having presented the general procedure to add dissipation in a system, we concentrate on the dissipative version of a model for cooperative behavior which, in its reversible version, has been extensively studied in \cite{Daw83} and, more recently, in \cite{gpy}. This model, characterized by a quartic self-potential and a mean-field interaction, shares many features of interacting FitzHugh-Nagumo models. We show that the thermodynamic limit can be studied in a suitable small-noise approximation, which reveals various peculiar features that we find interesting.
\bi
\item
A stable periodic orbit arises through a {\em homoclinic} bifurcation rather than a Hopf bifurcation; this is a {\em global} phenomenon, that is not detected by a local analysis.
\item
The stable periodic orbit may coexist with stable fixed points: the long-time behavior depends on the initial condition.
\item
The system can be {\em excited} by noise: by increasing the noise, stable fixed points may be de-stabilized, and trajectories are attracted by the periodic orbit which can be quite far from the fixed point.
\ei

We remark that the small noise approximations we obtain are quite similar to the class of Gaussian {\em nonlinear} processes originally studied in \cite{Sc1} and \cite{Sc2}, and more recently in \cite{touboul2011noise}. In particular, \cite{Sc1} and \cite{Sc2} give an interesting overview of oscillatory phenomena in these systems, and  of the role of the noise.

In Section \ref{prel} we describe a general class of dissipative models we consider. In Section \ref{propch} we derive the macroscopic limit of a specific class of models for cooperative behavior. The behavior of the macroscopic dynamics is studied in Section \ref{main}; in particular, we study in details the system in absence of noise and then we consider a suitable small noise approximation of the macroscopic equation, for which we prove noise-excitability. 
Most proofs are then postponed to Section \ref{proof} and to the Appendix.

\section{Systems driven by a random potential} \label{prel}

\subsection{Reference model}

Let ${\bf x}(t) = (x_i(t))_{i=1}^N$ denote the positions at time $t$  of $N$ particles moving in   $\R^d$. Particles are subject to a potential field, consisting of two parts: the external potential and  the interaction potential. The external potential is a smooth function $U: \R^d \ra \R$. Moreover, particles themselves generate a potential field, which is felt by all other particles. We adopt the mean-field assumption that particles contribute symmetrically to the field. More specifically, letting $\mathcal{M}_1(\R^d)$  be the space of Borel probabilities of $\R^d$ provided with the weak topology, and denoting by
\[
\rho_N({\bf x}(t)) := \frac{1}{N} \sum_{i=1}^N \d_{x_i(t)} \ \in \ \mathcal{M}_1(\R^d)
\]
the {\em empirical measure} of the particles at time $t$, we assume the potential at $x \in \R^d$ generated by the particles to be of the form
\be{basicV}
V_t(x) = \overline{V}(x, \rho_N({\bf x}(t)))
\ee
for some sufficiently regular function $\overline{V}: \R^d \times \mathcal{M}_1(\R^d) \ra \R$. Finally, for a given family $(w_i(t): t \geq 0 )_{i=1}^N$ of independent Brownian motions, we consider the evolution given by the following system of stochastic differential equations
\be{basicdyn}
dx_i(t) = - \nabla U(x_i(t)) dt - \nabla V_t(x_i(t)) dt + \sigma dw_i(t) \qquad (\mbox{for } i=1, \dots, N).
\ee

\subsection{Introducing dissipation and diffusion}

In the model above, the interaction potential in \eqref{basicV} is  a deterministic function of particle positions. When dissipative effects are present, as in a great variety of models in biology (see e.g. \cite{elowitz, garcia2004modeling, Schweitzer, touboul2011noise}) this scheme should be modified, by letting the potential $V_t$ evolving according to the following (Stochastic) Partial Differential Equation
\be{dissV}
dV_t(x) = - \a V_t(x) dt + D \Delta V_t(x) + d\overline{V}(x, \rho_N({\bf x}(t))) + \mbox{noise},
\ee
where $\a,D \geq 0$ tune dissipation and diffusion, respectively. The noise term, which could have several different forms, will be removed in the examples below. 

\br{rem:spin}
A similar construction can be applied to systems with discrete state space, e.g. spin systems for which $x_i \in \{-1,1\}$. In this case the Laplacian is omitted in the evolution \eqref{dissV} of $V_t$, and the evolution \eqref{basicdyn} is replaced by a Glauber dynamics for the time-dependent potential $\sum_i \left(U(x_i) + V_t(x_i) \right)$ at some inverse temperature $\b$. The special model corresponding to $U \equiv 0$ and
\[
\overline{V}(x, \rho_N({\bf x})) = - x \frac{1}{N} \sum_{i=1}^N x_i
\]
is studied in \cite{dfr}.

\er

\section{A model of cooperative behavior}
\label{propch}

We consider here $d=1$ and the following choice of external and interaction potential:
\[
U(x) = \frac{x^4}{4} - \frac{x^2}{2},
\]
and
\[
\overline{V}(x, \rho_N({\bf x})) = \frac{\theta}{2N} \sum_{i=1}^N (x-x_i)^2 \,,
\]
where $\theta$ is a positive parameter that represents the strength of the interaction between particles. The corresponding  reference (non-dissipative) model has been extensively studied in \cite{Daw83}. More recently, the same model has been applied to describe {\em systemic risk} in finance (see \cite{gpy}). It is known that in the limit as $N \ra +\infty$ the system exhibits a {\em phase transition}: for small $\theta$ the mean particle position is zero, while for large $\theta$ particles tend to cluster around the two minima of the double-well potential $U(x)$, giving rise to ``polarized'' equilibria.

To avoid complications, we assume the initial condition $V_0(x)$ to equation \eqref{dissV} is a polynomial  in $x$:
\[
V_0(x) := \sum_{k=0}^n a_k (0) x^k.
\]
Observing that, by Ito's rule
\[
 d\overline{V}(x, \rho_N({\bf x}(t))) = \theta x \, d\left( \frac{1}{N} \sum_{i=1}^N x_i(t) \right) - \frac{\theta}{N} \sum_{i=1}^N x_i(t) \, dx_i(t) + \frac{\s^2}{2} dt,
 \]
 equation \eqref{dissV} is solved by the polynomial
 \[
 V_t(x) := \sum_{k=0}^n a_k(t) x^k
 \]
 whose coefficients satisfy the system of equations
 \be{coeff}
 \left\{
 \begin{array}{rcll}
 da_k(t) & = & - \a a_k(t) dt & \mbox{for } k= n,n-1 \\
 da_k(t) & = & - \a a_k(t) dt + D(k+2)(k+1)a_{k+2}(t) dt& \mbox{for } 2 \leq k \leq n-2 \\
 da_1(t) & = & - \a a_1(t) dt + 6Da_{3}(t) dt + \theta d\left( \frac{1}{N} \sum_{i=1}^N x_i(t) \right) .& 
 \end{array}
 \right.
 \ee
 Note that $a_0(t)$ is not needed for \eqref{basicdyn}. By \eqref{coeff}, it follows that $a_k(t) \ra 0$ as $t \ra +\infty$, for all $k \geq 2$. Since we are interested in a steady state regime, we may assume $V_t(x)$ is linear in $x$ for all $t$. Renaming $\mu(t) := a_1(t)$, we obtain the following form for \eqref{basicdyn} coupled with \eqref{dissV}:
\be{micro}
\begin{array}{rcl}
dx_i(t) &=& \left( -x_i^3(t) + x_i(t) \right) dt - \mu(t) dt + \sigma  dw_i(t) \qquad (\mbox{for } i=1, \dots, N)\\ & &  \\
d\mu (t) &=& -\alpha \mu(t) dt - \theta dm^{(N)}(t) \,,
\end{array}
\ee
where we set $ m^{(N)}(t)  :=  \frac{1}{N} \sum_{i=1}^N x_i(t) $. Note that \eqref{micro} can be rewritten as
\be{micro1}
\begin{array}{rcl}
dx_i(t) &=& \left( -x_i^3(t) + x_i(t) \right) dt - \mu(t) dt + \sigma  dw_i(t) \qquad (\mbox{for } i=1, \dots, N)\\ & &  \\
d\mu (t) &=& -(\alpha - \theta) \mu(t) dt - \frac{\theta}{N} \sum_{i=1}^N \left( -x_i^3(t) + x_i(t) \right) dt - \frac{\theta \s}{N} \sum_{i=1}^N dw_i(t) \,.
\end{array}
\ee
We remark that, although the drift terms in \eqref{micro1} are not uniformly Lipschitz, strong existence and uniqueness can be established, for example by using the classical Hasminskii's Test: by defining, for instance
\[
V({\bf x}, \mu) := \frac{1}{N} \sum_{i=1}^N \left[ \frac{x_i^4}{4} + \frac{x_i^2}{2} \right] + \frac{a}{2} \mu^2
\]
for $a>0$ sufficiently small, one obtains an inequality of the form
\[
\mathcal{L}V({\bf x},\mu) \leq K\left[1+V({\bf x},\mu) \right]
\]
for some $K>0$, where $\mathcal{L}$ is the infinitesimal generator of the diffusion \eqref{micro1}. This inequality implies existence and uniqueness of strong solutions (see e.g. \cite{mk}).

\subsection{Propagation of chaos}

In the limit as $N \ra +\infty$, the system of equations \eqref{micro1} naturally suggest the following limiting process, which describes the behavior of a single component in the limit:
\be{macro}
\begin{split}
dx(t) & =  \left( -x^3(t) + x(t) \right) dt - \mu(t) dt + \s dB(t) \\ 
\frac{d\mu}{dt}(t) & =  - (\a - \theta) \mu(t) + \theta \E \left[x^3(t) - x(t) \right],
\end{split}
\ee
where $B(\cdot)$ is a standard Brownian motion. Note that this equation has a {\em non-local} structure, due to the appearance of the law of $x(t)$ in the term $ \E \left[x^3(t) - x(t) \right]$ of the right hand side. 
\bt{th:ex}
For each $\mu_0 \in \R$ and each real random variable $\xi$, with finite third moment and independent of the Brownian motion $B(\cdot)$, equation \eqref{macro} with initial condition $x(0) = \xi$, $\mu(0) = \mu_0$, has a unique strong solution.
\et
\bt{th:prop-ch}
Let $(x_i^{(N)}(\cdot), \mu^{(N)}(\cdot))_{i=1}^N$ be the solution of \eqref{micro1} with an initial condition satisfying the following conditions:
\bi
\item[(a)]
$(x_i^{(N)}(0))_{i=1}^N$ are independent and identically distributed with a given law $\l$ having finite third moment; moreover, they are independent of the Brownian motions $(w_i(\cdot))_{i=1}^N$.
\item[(b)]
$\mu^{(N)}(0) = \mu_0 \in \R$.
\ei
Then, as $N \ra +\infty$ and for each fixed $k \geq 1$, the vector-valued process $$(x_1^{(N)}(\cdot), x_2^{(N)}(\cdot), \ldots, x_k^{(N)}(\cdot))$$ restricted to a given time interval $[0,T]$ converges in law to $(y_1(\cdot), y_2(\cdot), \ldots, y_k(\cdot))$, where $y_1(\cdot), y_2(\cdot), \ldots, y_k(\cdot)$ are independent copies of the solution of \eqref{macro} with initial condition $y_j(0) \sim \l$, $\mu(0) = \mu_0$.

\et

The proofs of Theorems \ref{th:ex} and \ref{th:prop-ch} are rather standard; for completeness, we sketch them in the Appendix. Note that the only difficulty is to deal with the non-global Lipschitz property of the drift. In the microscopic model  this problem is overcome by stopping the process at the boundary of a compact set and controlling this stopping time via Lyapunov methods. This approach is not directly applicable at the macroscopic level: stopping the process affects the drift globally, due to the  non-locality of the evolution.

\section{Main results}

\label{main}

Our main goal is to study the long-time evolution of solutions to equation \eqref{macro} and of the corresponding Fokker-Planck equation:
\be{FP}
\begin{array}{rcl}
\partial_t q_t (x) &=& \dfrac{\sigma^2}{2} \partial^2_{xx} q_t (x) - \partial_x \left\{ \left[ -x^3 + x - \mu(t) \right] q_t (x)\right\} \\
\dfrac{d\mu}{dt}(t) &=& - (\alpha - \theta) \mu(t) - \theta \left\langle -x^3 + x, q_t \right\rangle \,,
\end{array}
\ee
where the notation $\langle f,q \rangle := \int f(x) q(x) dx$ is used.
The regularizing effect of the second derivative guarantees that, for $t>0$, the law of $x(t)$ has a density $q_t(x)$ which solves \eqref{FP}.

\subsection{Stationary solutions}

The study of equilibria for the macroscopic dynamics turns out to be essentially trivial.
\bp{prop:equilibrium}
Equation \eqref{FP} admits a unique stationary solution given by
\begin{equation}\label{StatSol:MacroDyn}
q_* (x) = Z^{-1}_* \cdot \exp \left\{ \frac{1}{\sigma^2} \left[ - \frac{x^4}{2} + x^2 \right] \right\} \,, \quad \mu_* = 0
\end{equation}
where $Z_*$ is a normalizing factor.

\ep
We remark that the scenario is quite different than for the reference model with the same potential but without dissipation, where multiple equilibria arise for large $\theta$.

\subsection{The noiseless dynamics}

Letting $\s = 0$ in \eqref{macro}, we obtain a deterministic dynamics described by the ODE
\be{macrodet}
\begin{split}
\frac{dx}{dt}(t) & =   -x^3(t) + x(t)  - \mu(t)  \\ 
\frac{d\mu}{dt}(t) & =  - (\a - \theta) \mu(t) + \theta \left[x^3(t) - x(t) \right].
\end{split}
\ee
The dependence of the attractors for \eqref{macrodet} on the parameters is quite nontrivial, despite of  the fact that \eqref{macrodet} has always $(0,0)$, $(-1,0)$ and $(1,0)$ as unique equilibrium points in the plane $(x,\mu)$. Moreover, for all values of the parameters $\a > 0$ and $\theta \geq 0$, the equilibrium $(0,0)$ is a saddle point. We denote by $W$ the {\em stable manifold} of the origin, i.e. the set of initial conditions whose corresponding solution converges to the origin as $t \ra +\infty$.

\bt{th:set}
Fix $\a>0$. Then there is $\theta_1 = \theta_1(\a) > 0$ such that the following properties hold.
\bi
\item[(a)]
For $\theta < \theta_1$, $W$ is an unbounded, open curve, that separates the basins of attraction of $(-1,0)$ and $(1,0)$ respectively.
\item[(b)]
For $\theta = \theta_1$, $W$ closes into an eight-shaped curve, comprised by two cycles, symmetric with respect to the 
origin, surrounding $(-1,0)$ and $(1,0)$ respectively (homoclinic bifurcation).
\item[(c)]
For $\theta_1 < \theta < \a+2$, $W$ spirals around two unstable periodic solutions of \eqref{macrodet}, surrounding $(-1,0)$ and $(1,0)$ respectively. These cycles separate the basins of attraction of $(-1,0)$ and $(1,0)$ from the basin of attraction of a unique stable periodic solution of \eqref{macrodet}, which surrounds both $(-1,0)$ and $(1,0)$.
\item[(d)]
At $\theta = \a+2$, the two unstable cycles collapse in $(-1,0)$ and $(1,0)$ through a Hopf bifurcation. For $\theta \geq \a+2$, $(-1,0)$ and $(1,0)$ become unstable: \eqref{macrodet} admits a unique periodic solution, whose basin of attraction is $\R^2 \setminus \left[\{(\pm 1,0)\} \cup W \right]$.
\ei

\et

\subsection{Excitability by noise: Gaussian approximation}

\begin{figure}[h!]
\centering%
\subfigure{\includegraphics[width=0.47\textwidth]{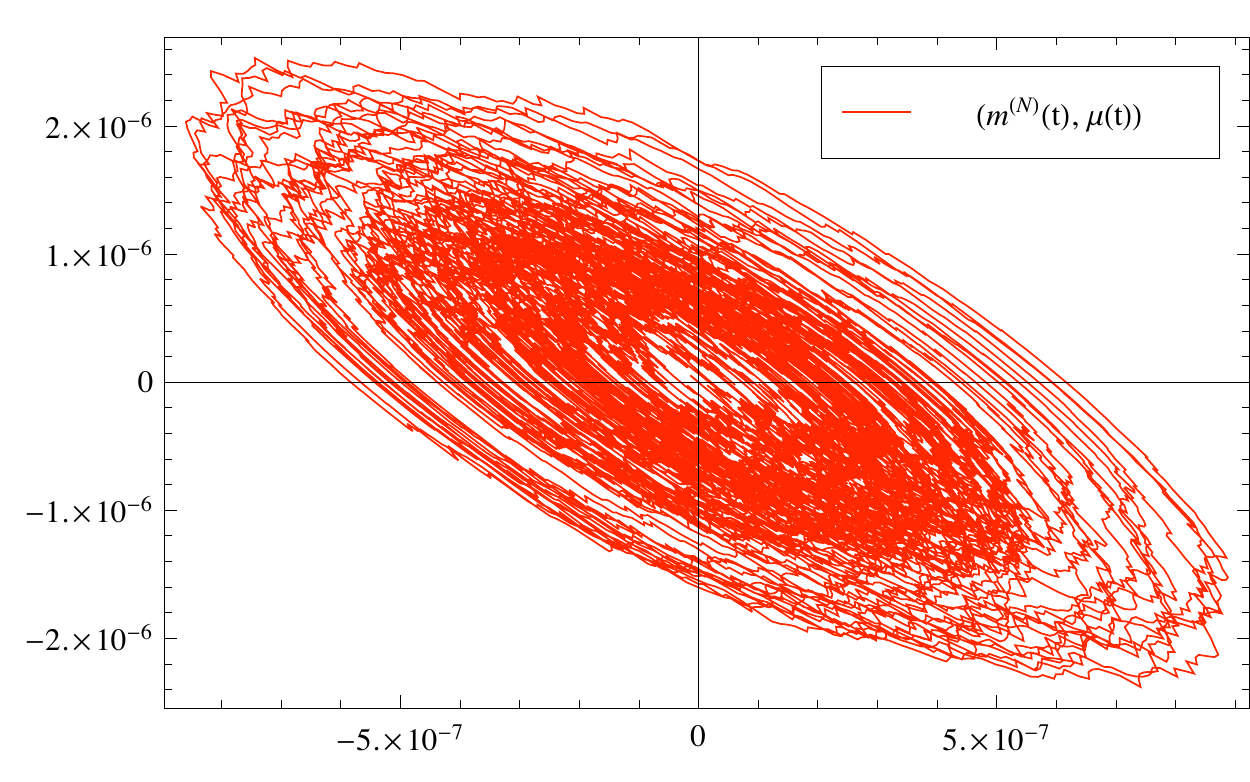}} 
\subfigure{\includegraphics[width=0.45\textwidth]{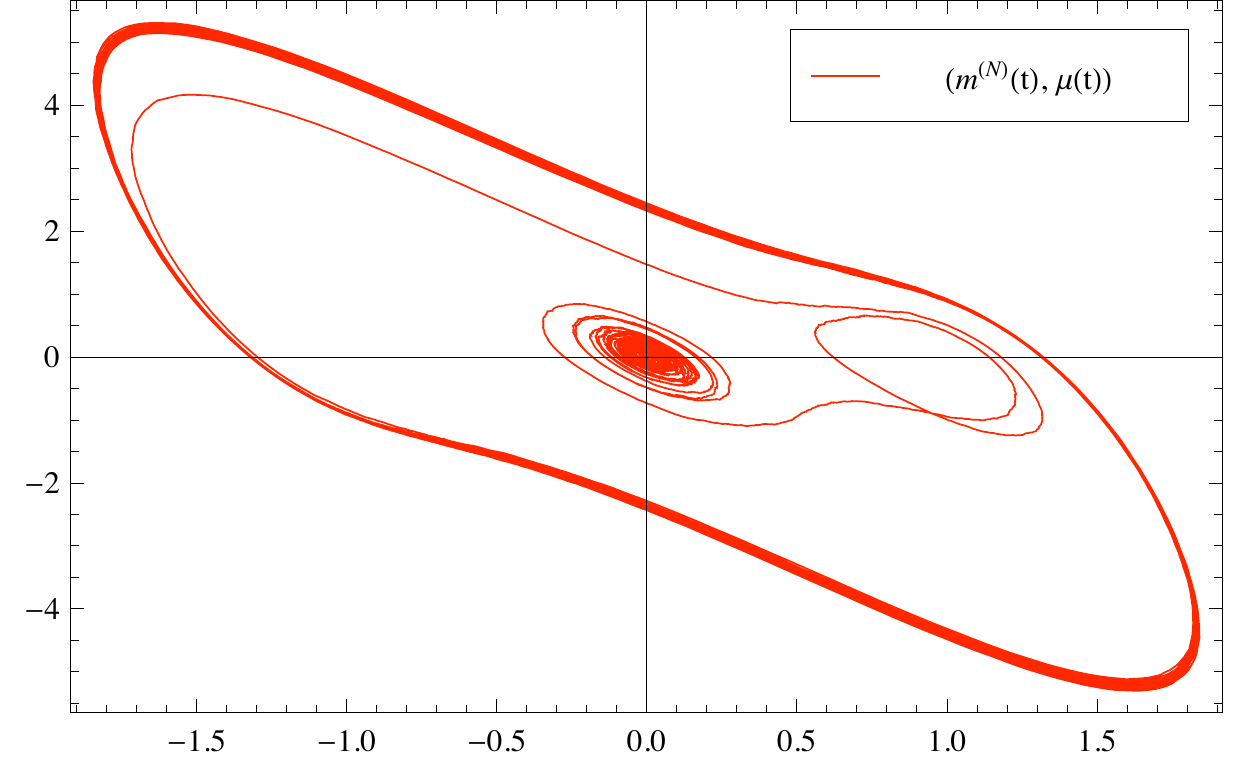}}
\caption{
The figures illustrate the phenomenon of \emph{excitation by noise}, that is how noise (of suitable intensity) helps/favors the emergence of a stable periodic behavior. Simulations of the interacting particle system \eqref{micro} have been run for $\theta$ slightly below the threshold value $\alpha+2$. Recall that this range of the parameter corresponds to a coexistence phase for the deterministic system \eqref{macrodet}, in which the locally stable equilibria $(\pm 1, 0)$ coexist with a macroscopic stable limit cycle. We simulate the trajectory $t \mapsto (m^{(N)}(t),\mu(t))$ for different noise sizes. 
Both simulations are performed for the same number of iterations, with the same time-step size, and starting from $(m^{(N)}(0),\mu(0))=(0,0)$ with particles initially placed half in $(+1,0)$ and half in $(-1,0)$.  
Only the noise size is changing. On the left panel, when $\sigma$ is small, the curve $(m^{(N)}(t),\mu(t))$ stays close to the initial point for all times; indeed, it oscillates with maximum amplitude of order of $10^{-6}$. On the right panel, we increase the noise size and $(m^{(N)}(t),\mu(t))$ exhibits a periodic behavior. 
}
\label{fig:Small:Big:Noise}
\end{figure}

The study of the effects of the noise ($\s>0$) in \eqref{macro} appears to be of a higher level of difficulty compared to the analysis of the deterministic system. 
Simulations of the interacting particle system dynamics \eqref{micro}, indicate that the following picture arises. 
\bi
\item
The structure of the attractors is preserved for small $\s$, except that the fixed points $(\pm 1,0)$ are replaced by the invariant law \eqref{StatSol:MacroDyn}. For $\s$ sufficiently large, the noise breaks the organized structure, so that, in particular, there is no periodic solution.
\item
For intermediate values of $\s$ the behavior depends on the parameters. The most interesting case is for $\theta_1 < \theta < \a+2$. By increasing $\s$, the basin of attraction of the invariant law \eqref{StatSol:MacroDyn} shrinks, and this invariant law becomes unstable at a critical value $\s = \s_c$. A unique stable periodic solution survives for $\s < \s'_c$, with $\s'_c > \s_c$. Thus, an initial condition that, for small $\s$, would be in the basin of attraction of \eqref{StatSol:MacroDyn}, gets {\em excited} by the noise, and is attracted by a periodic orbit. See Figure~\ref{fig:Small:Big:Noise}.
\ei

We have no proof of the above statements, but strong evidences is obtained through numerical simulations. Surprisingly, even the analysis of the linear stability of \eqref{StatSol:MacroDyn} appears to be quite hard. One possible approach consists in considering the infinite dimensional system of equations for the moments of \eqref{macro}. By Ito's rule we obtain, for $k \geq 1$:
\[
dx^k(t) 
= \left[ -k x^{k+2}(t) + k x^k(t) + \frac{\sigma^2}{2} k (k-1) x^{k-2}(t) -k \mu(t) x^{k-1}(t) \right] dt + \sigma k x^{k-1} (t) dB(t) \,.
\]
Setting $m_k(t) := \E(x^k(t))$, taking the expectation in the previous equation we get
\be{moments}
\begin{split}
\dfrac{dm_k}{dt}(t) & = -k m_{k+2}(t) + k m_k(t) + \dfrac{\sigma^2}{2} k (k-1) m_{k-2}(t) -k \mu(t) m_{k-1}(t)  \\
\frac{d\mu}{dt}(t) & = - (\alpha-\theta) \mu(t) - \theta \left( -m_3(t) + m_1(t)\right) = -\a \mu(t) - \theta \dfrac{dm_1(t)}{dt}.
\end{split}
\ee
This system of equations appears to be as hard as \eqref{FP}. However, we use it as a guide to construct a 
Gaussian approximation of \eqref{macro}. To be precise, denote by $(x_{\s}(t), \mu_{\s}(t))$ the solution of \eqref{macro}, for a given deterministic initial condition $x_{\s}(0) = x, \mu_{\s}(0) = \mu$ (random initial conditions could be considered as well). We construct, on the same probability space, a process $(y_{\s}(t), \nu_{\s}(t))$ satisfying the following conditions:
\bi
\item[{\bf (a)}]
$y_{\s}(\cdot)$ is a Gauss-Markov process, $\nu_{\s}(\cdot)$ is a deterministic process, $y_{\s}(0) = x$, $\nu_{\s}(0) = \mu$.
\item[{\bf (b)}]
The moments of $y_{\s}(\cdot)$  and $\nu_{\s}(\cdot)$ satisfy \eqref{moments} for $k=1,2$.
\item[{\bf (c)}]
For every $T>0$ there is a constant $C_T$ such that for all $\s>0$
\be{normappr}
\E \left[ \sup_{t \in [0,T]} |x_{\s}(t) - y_{\s}(t)| \right] +  \sup_{t \in [0,T]} |\mu_{\s}(t) - \nu_{\s}(t) | \leq C_T \s^2.
\ee
\ei
Note that \eqref{normappr}, which expresses the fact that $y_{\s}$ is a small noise approximation of $x_{\s}$, is not enough to identify the law of $y_{\s}$, since this condition in insensitive to corrections of order $\s^2$. These corrections  are determined by condition {\bf (b)}; indeed, the equations \eqref{moments} for $k=1,2$ form a closed system of two equations, since for Gaussian random variables, all moments can be expressed in terms of the first two.

We now define the Gaussian process $y_{\s}$. First, consider the following system of ODEs:
\be{gaussapprox}
\begin{split}
\frac{dm_{\s}}{dt} & = - m_{\s}^3 + m_{\s}- \nu_{\s} - 3 \s^2 m_{\s} V_{\s} \\
\frac{d\nu_{\s}}{dt} & = - \a \nu_{\s} - \theta \frac{dm_{\s}}{dt} \\
\frac{dV_{\s}}{dt} &= 1 + 2(1-3m_{\s}^2)V_{\s} - 6 \s^2 V_{\s}^2 \\
m_{\s}(0) &= x \ \ \ \nu_{\s}(0) = \mu \ \ \ V_{\s}(0) = 0.
\end{split}
\ee
It is easily seen that solutions of \eqref{gaussapprox} are bounded, which implies global existence and uniqueness. We can therefore define the following Gauss-Markov, centered process, as the unique solution of the linear SDE:
\begin{equation}\label{linearsde}
\begin{split}
dz_{\s}(t) & = \left[1-3m_{\s}^2(t) - 3 \s^2 V_{\s}(t) \right] z_{\s}(t) dt + dB(t) \\
z_{\s}(0) & = 0,
\end{split}
\end{equation}
where $B(t)$ is the same Brownian motion driving \eqref{macro}. By computing $dz_{\s}^2$ one checks that
\[
Var(z_{\s}(t)) = V_{\s}(t).
\]
Finally, let
\be{gaussY}
y_{\s}(t) := m_{\s}(t) + \s z_{\s}(t).
\ee

\bp{prop:gaussapprox}
Properties {\bf (a)}, {\bf (b)} and {\bf (c)} hold for the process $y_{\s}$ defined in \eqref{gaussY}.
\ep

\begin{figure}[h!]
\centering%
\subfigure{\includegraphics[width=0.45\textwidth]{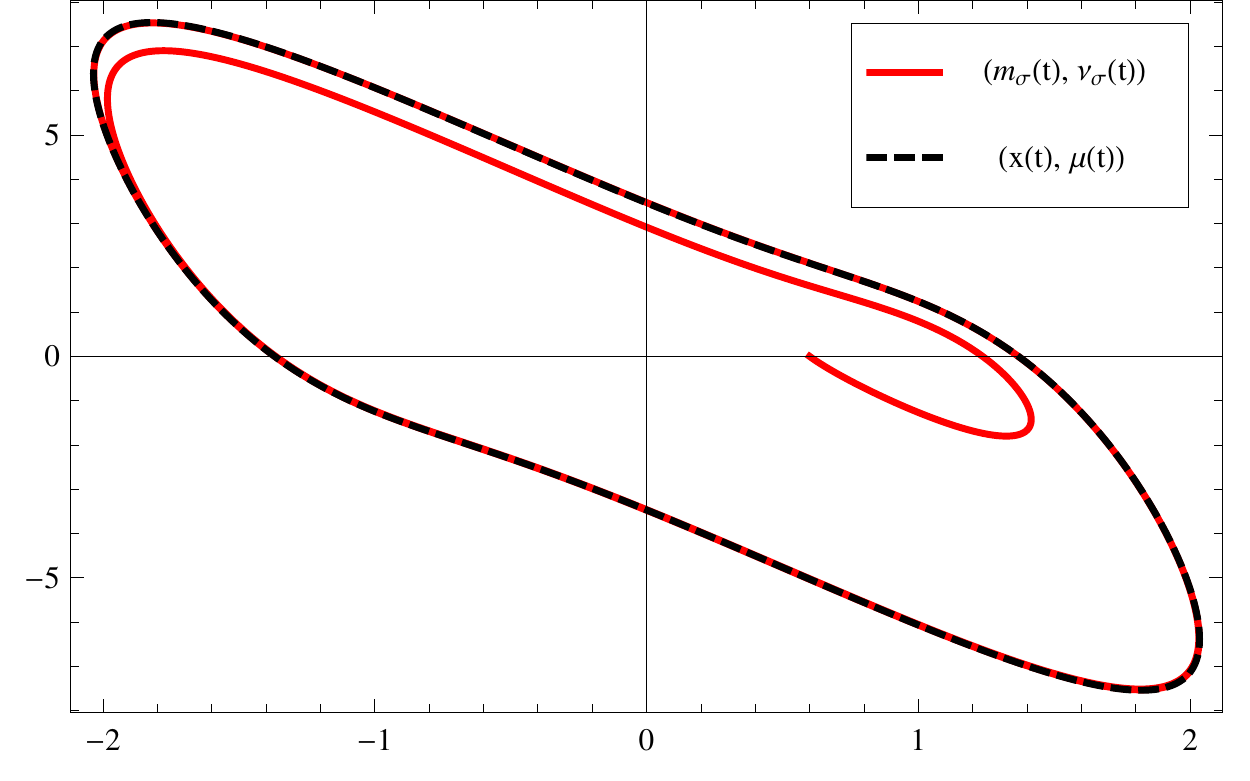}} 
\subfigure{\includegraphics[width=0.45\textwidth]{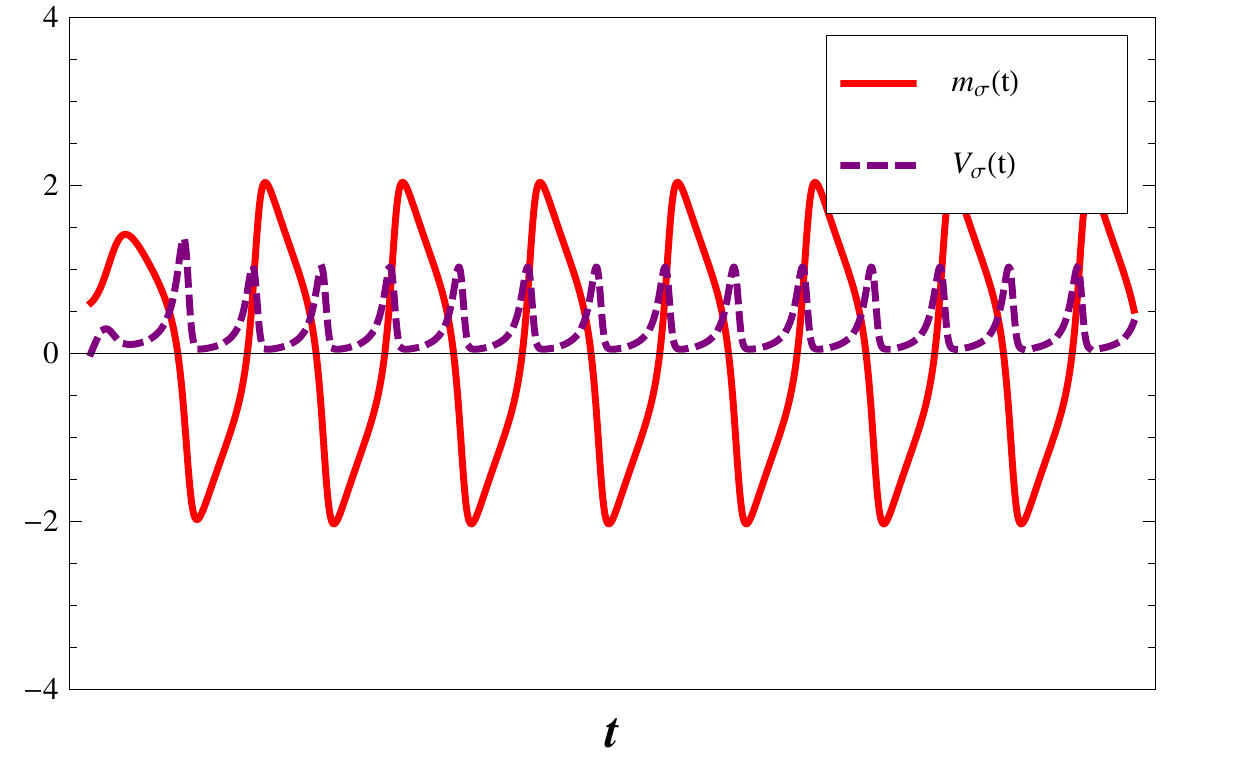}}
\caption{
Regime $\sigma\ll 1$ and $\theta>\alpha+2$. In the left panel the \mbox{$(V_\sigma = 0)$--section} of the typical motion of $(m_{\sigma},\nu_{\sigma})$ is shown.  The trajectory $t \mapsto (m_{\sigma}(t),\nu_{\sigma}(t))$ (solid red line) approaches and remains close to the limit cycle attractor for the dynamics \eqref{macrodet} (dashed black line). A heuristic explanation of the phenomenon can be inferred from the right panel, where the oscillatory behavior of $m_{\sigma}$ and $V_{\sigma}$ is illustrated. Since $V_{\sigma}$ is of order one for all times, the terms proportional to  $\sigma^2V_{\sigma}$ in \eqref{gaussapprox} are small when $\sigma\ll 1$. Thus, the trajectories of \eqref{gaussapprox} are close to those of the deterministic system \eqref{macrodet}. 
}
\label{fig:Mean:Var:Osc}
\end{figure}

Having obtained a small noise Gaussian approximation for \eqref{macro}, we study the approximated process. Being Gaussian, it is enough to study the evolution of mean and variance. Since $\E(y_{\s}(t)) = m_{\s}(t)$ and $Var(y_{\s}(t)) = \s^2 V_{\s}(t)$, we are left to the analysis of the three dimensional system \eqref{gaussapprox}. In the following analysis we omit the subscript $\s$ in $m_{\s}, \nu_{\s}, V_{\s}$. In the three dimensional space $\{(m,\nu,V): \, m,\nu \in \R, \, V \geq 0\}$, equations \eqref{gaussapprox} admit the following equilibrium points:
\[
s_1=\left(  \frac{\sqrt{1+\sqrt{1-3 \sigma^2}}}{\sqrt{2}},\: 0,\: \frac{1-\sqrt{1-3\sigma^2}}{6\s^2} \right), \quad 
s_2=\left( - \frac{\sqrt{1+\sqrt{1-3 \sigma^2}}}{\sqrt{2}},\: 0,\: \frac{1-\sqrt{1-3\sigma^2}}{6\s^2} \right),
\]
\[
s_3  =\left( - \frac{\sqrt{1-\sqrt{1-3 \sigma^2}}}{\sqrt{2}},\: 0,\: \frac{1+\sqrt{1-3\sigma^2}}{6\s^2} \right), \quad
s_4 = \left( \frac{\sqrt{1-\sqrt{1-3 \sigma^2}}}{\sqrt{2}},\: 0,\: \frac{1+\sqrt{1-3\sigma^2}}{6\s^2} \right),
\]
for $\s^2 \leq 1/3$, and
\[ 
s_5 = \left( 0,\: 0,\: \frac{1+\sqrt{1+6\sigma^2}}{6\s^2}\right)
\]
for all $\s>0$. Note that, as $\s \downarrow 0$, the equilibria $s_1$ and $s_2$ converge to $(\pm 1,0,3/2)$, so they correspond to the equilibria $(\pm 1,0)$ of the deterministic system \eqref{macrodet}; the equilibria $s_3$, $s_4$, and $s_5$ converge to $(0,0,+\infty)$. In the following result we give a {\em local} version of the property of excitability by noise: under suitable conditions on the parameters, by increasing $\s$ the equilibria $s_1$ and $s_2$ lose stability.
%
\bt{th:excitable}
Assume $\theta < \a +2$, with $\a + 2 - \theta$ sufficiently small, and assume $\a<10$. Then there exists $\s_c>0$ such that $s_1$ and $s_2$ are linearly stable for $\s<\s_c$, but linear stability is lost as $\s$ crosses $\s_c$.
\et
For $\theta > \alpha + 2$ and $\sigma \ll 1$, we expect system \eqref{gaussapprox} to have an unique attracting periodic solution, as the deterministic system \eqref{macrodet} has. This behavior is illustrated by simulations in Figure~\ref{fig:Mean:Var:Osc}.

\begin{figure}[h!]
\centering%
\subfigure{\includegraphics[width=0.45\textwidth]{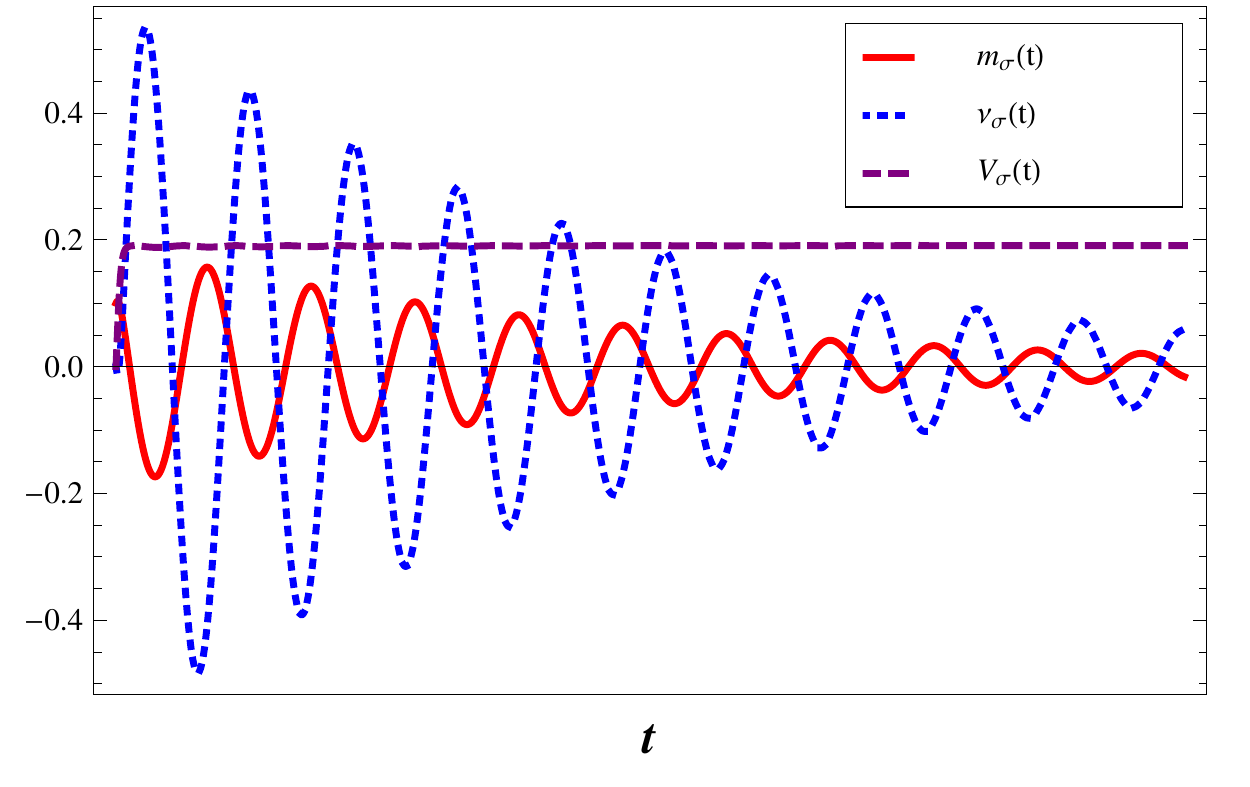}} 
\subfigure{\includegraphics[width=0.45\textwidth]{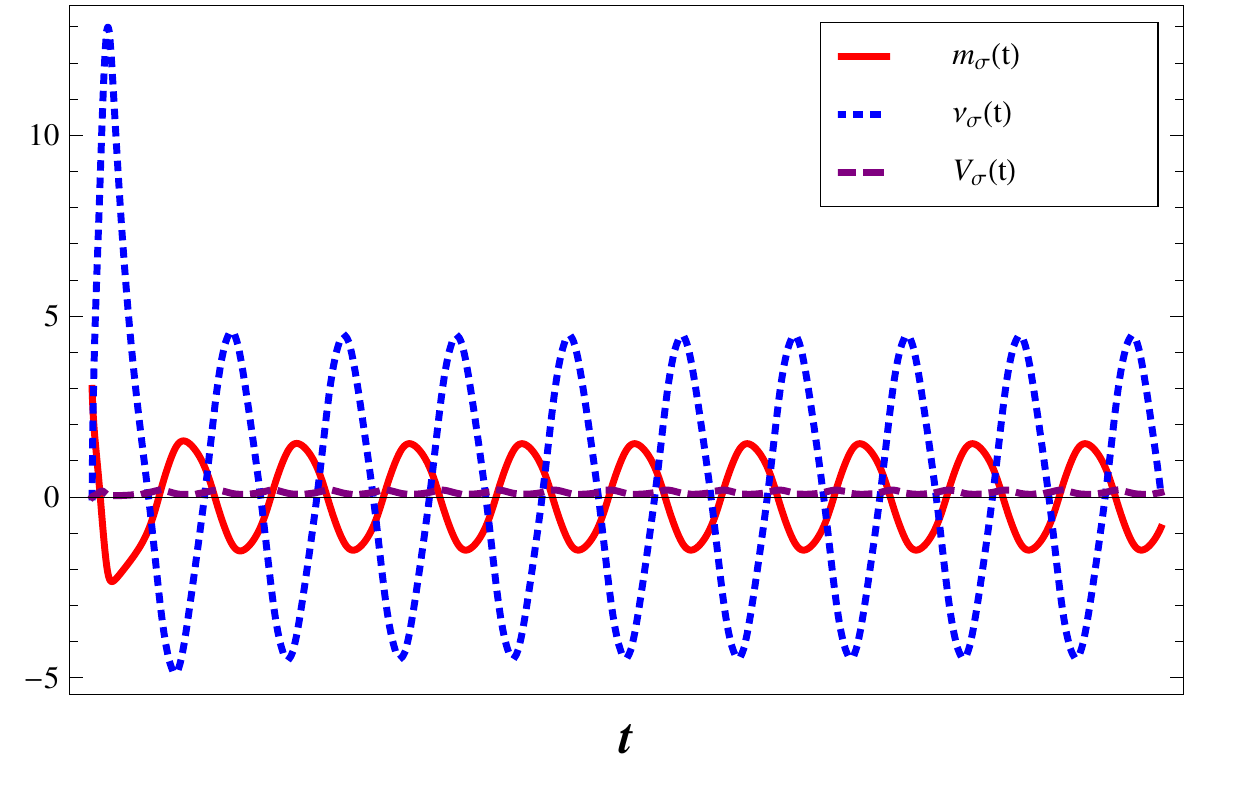}}
\caption{
Simulations for the system \eqref{gaussapprox} in the case when $\theta>\alpha+2$ and $\sigma>\sqrt{\frac{2}{3}(\alpha-\theta)(\alpha-\theta-1)}$. In this regime the limiting behavior of the trajectories depends on the initial value of the magnetization. We run simulations with starting point $V_{\sigma}(0)=\nu_{\sigma}(0)=0$ and various initial conditions for $m_{\sigma}(0)$. 
For $\left |m_{\sigma}(0)\right| \ll 1$ the oscillatory behavior disappears. Indeed, $(m_{\sigma},\nu_{\sigma},V_{\sigma})$ converges to $s_5=(0,0,\frac{1+\sqrt{1+6\sigma^2}}{6 \sigma^2})$, which is the only fixed point for the dynamics \eqref{gaussapprox} in this phase and is linearly stable (left panel). On the contrary, when $\left|m_{\sigma}(0)\right|$ is large, the periodicity persists also for $\sigma\gg 1$ (right panel).
}
\label{fig:Osc:Sigma:Grande}
\end{figure}

The behavior of \eqref{gaussapprox} for large $\s$ could be studied as well, although it may not reflect any of the features of \eqref{macro}. For $\theta > \alpha + 2$ one shows, in particular, that the equilibrium $s_5$ becomes linearly stable as $\sigma>\sqrt{\frac{2}{3}(\alpha-\theta)(\alpha-\theta-1)}$, but it is not a global attractor, since a locally stable periodic orbit survives. This is illustrated in Figure \ref{fig:Osc:Sigma:Grande}.

\section{Proofs}

\label{proof}

\subsection{Proof of Proposition \ref{prop:equilibrium}}

In view of dynamics \eqref{FP}, an equilibrium probability density for our system must satisfy the system of equations
\begin{subequations}\label{StatEq}
\begin{equation}\label{StatEq:q}
\frac{\sigma^2}{2} \partial^2_{xx} q_* (x) - \partial_x \left\{ \left[ -x^3 + x - \mu_* \right] q_* (x)\right\} = 0
\end{equation}
\begin{equation}\label{StatEq:mu}
- (\alpha - \theta) \mu_* - \theta  \left\langle -x^3 + x, q_* \right\rangle = 0 \,. 
\end{equation}
\end{subequations}
Solving \eqref{StatEq:q} we obtain
\[
q_* (x) = Z_*^{-1} \cdot \exp \left\{ \frac{1}{\sigma^2} \left[ - \frac{x^4}{2} + x^2 + 2 \theta \mu_* x \right]\right\}
\]
with
\[
Z_* = \int_{\mathbb{R}} \exp \left\{ \frac{1}{\sigma^2} \left[ - \frac{x^4}{2} + x^2 + 2 \theta \mu_* x \right]\right\} dx \,.
\]
This result together with equation \eqref{StatEq:mu} imply $\mu_*=0$ and the proof is complete.

\subsection{Proof of Theorem \ref{th:set}}

The scenario depicted in the statement is due to the occurrence of both a homoclinic and a Hopf bifurcation. To ease the readability of the proof, we believe it is worth first to explain what is going on and then to give the technical details. \\
Roughly speaking, the theorem states there exist three possible phases for system \eqref{macrodet}:
\begin{itemize}
\item
\emph{Fixed points phase}. For $\theta<\theta_1$ the only stable attractors are $(\pm 1,0)$. 
\item
\emph{Coexistence phase}. For $\theta=\theta_1$ the system has an eight-shaped, symmetric with respect to the origin, homoclinic orbit surrounding $(\pm 1,0)$. By increasing the parameter $\theta$ from $\theta_1$, this separatrix splits into an outer stable limit cycle, which  contains both the equilibria, and two inner unstable periodic orbits around $(1,0)$ and $(-1,0)$, respectively. In this phase $(\pm 1, 0)$ are linearly stable. Therefore, two locally stable fixed points coexist with a stable limit cycle. 
\item
\emph{Periodic orbit phase}. For $\theta=\alpha+2$ a Hopf bifurcation occurs: the inner unstable limit cycles disappear collapsing at $(\pm 1,0)$. At the same time the equilibria lose their stability and, thus, the external stable limit cycle remains the only stable attractor for $\theta\geq\alpha+2$.
\end{itemize}
The key point of the proof is the particular structure of the vector field generated by \eqref{macrodet}. 
For convenience, let us denote such a vector field by 
\begin{equation}\label{vector:field}
V_{\alpha,\theta}(x,\mu) := 
\begin{pmatrix}
-x^3+x-\mu \\ 
-(\alpha-\theta)\mu +\theta (x^3-x)
\end{pmatrix}.
\end{equation} 
Observe that \eqref{vector:field} defines a \emph{one-parameter family of negatively rotated vector fields} (with respect to $\theta$, for fixed $\alpha$); that is, the following are satisfied 
\begin{itemize}
\item 
the critical points of \eqref{vector:field} are isolated; 
\item
at ordinary points, as the parameter $\theta$ increases, all the field vectors $V_{\alpha,\theta}(x,\mu)$ rotate clockwise.
\end{itemize} 
Additionally, \eqref{vector:field} is \emph{semi-complete}; in other words, the field rotates of an angle $\pi$ as the parameter $\theta$ varies over $\mathbb{R}$. \\ 
For dynamical systems that depend on this specific way on a parameter, many results concerning bifurcations, stability and global behavior of limit cycles and separatrix curves are known. For the sake of completeness and readability, we collect here the properties satisfied by \eqref{vector:field} that are crucial in the sequel, referring to \cite[Chapter 4]{Per01} for precise and general statements.  
\begin{enumerate}
\item \label{P1} 
Limit cycles expand or contract monotonically as the parameter $\theta$ varies in a fixed sense.
\item \label{P2} 
A limit cycle is generated/absorbed either by a critical point or by a separatrix of  \eqref{vector:field}.  
\item \label{P3}
Cycles of distinct fields do not intersect.\\
\end{enumerate}

Now let us go into the details of the proof.

\paragraph{Lyapunov function} Let $\alpha, \theta > 0$. Consider the function
\[
W(x,\mu) = \frac{x^2}{2} + \frac{(\theta x + \mu)^2}{2 \alpha \theta}
\]
and observe that the total derivative
\[
\frac{d}{dt} W(x,\mu) = -x^4 + (1+\theta) x^2 - \frac{1}{\theta} (\theta x + \mu)^2
\]
is negative for every $\mu$ and sufficiently large $x$. Thus, there exists a stable domain for the flux of \eqref{macrodet} and, in particular, the trajectories can not escape to infinity as $t \to +\infty$. 

\paragraph{Local analysis of equilibria and Hopf bifurcation} As already mentioned, the vector field \eqref{vector:field} admits three fixed points in the phase plane $(x,\mu)$: $(0,0)$ and $(\pm 1,0)$. The origin is a saddle point for all values of the parameters; while, the equilibria $(\pm1,0)$ are linearly stable for $\theta < \alpha + 2$ and their local stability is lost for $\theta > \alpha + 2$. At the critical point $\theta = \alpha + 2$ a Hopf bifurcation occurs: by decreasing the parameter $\theta$ from $\alpha + 2$ two symmetric unstable limit cycles are generated from $(\pm1,0)$. This proves the first assertion in statement (d).

\paragraph{Separatrix formation and stability} Properties \ref{P1} and \ref{P2} allow us to explain the appearance of a separatrix curve, whose breakdown causes a homoclinic bifurcation at $\theta=\theta_1$. Indeed, while decreasing $\theta$ from $\alpha+2$, the cycles arisen from $(\pm 1,0)$ expand until they join each other at the origin forming a homoclinic eight-shaped separatrix graphic $\Gamma_0$ at  $\theta=\theta_1$. 
Notice that, for $\theta=\theta_1$, $\Gamma_0 = W$. This proves statement (b).


Furthermore the exterior of the separatrix is stable. In fact, the existence of a stable domain for the flux of \eqref{macrodet} implies that all trajectories in an outer neighborhood of $\Gamma_0$ approach $\Gamma_0$ itself as $t \to +\infty$. 

\paragraph{Description of the phases} We are left to prove statements (a), (c) and part of (d). Hence,
\begin{itemize}
\item[(c)]
When the separatrix graphic $\Gamma_0$ splits increasing $\theta$ from $\theta_1$, it generates a periodic orbit surrounding both $(\pm 1,0)$ on its exterior and two smaller limit cycles, around $(-1,0)$ and $(1,0)$ respectively, on its interior. \\
The inner cycles are unstable (due to the subcritical Hopf bifurcation at $\theta=\alpha+2$) and represent the boundaries of the basins of attraction of $(\pm 1, 0)$. Moreover, the external periodic orbit inherits the stability of the exterior of $\Gamma_0$ and so it is stable. See Theorem 2 in \cite[Section 4.6]{Per01} for more details. 
 
\item[(a)]
It suffices to prove that in this phase the dynamical system \eqref{macrodet} does not admit a periodic solution. Indeed, the non-existence of cycles together with the existence of a stable domain for the flux of \eqref{macrodet} 
guarantee that every trajectory must converge to an equilibrium as $t \to +\infty$. The fact that the basins of attraction of $(\pm 1,0)$ are separated by $W$ is a consequence of the Stable Manifold theorem (see \cite[Section 2.7]{Per01}). \\
Thus, it remains to show that there is no limit cycle for $\theta < \theta_1$. From properties \ref{P1} and \ref{P2} it follows that, as $\theta$ increases from $\theta_1$ to infinity, the outer stable limit cycle expands and its motion covers the whole region external to $\Gamma_0$. Similarly, the two unstable cycles contract from the graphic $\Gamma_0$ and terminate at the critical points $(\pm1,0)$. As a consequence, for $\theta >\theta_1$ the entire phase space is covered by expanding or contracting limit cycles. Now, by using property \ref{P3}, we can deduce that  no periodic trajectory may exist  for $\theta < \theta_1$. In fact, such an orbit would intersect some of the cycles present when  $\theta>\theta_1$ that is not possible. 
\item[(d)] The statement easily follows by combining the impossibility of escaping trajectories, 
the presence of the stable manifold of the origin, the loss of stability of the fixed points $(\pm 1, 0)$ and the existence of a stable limit cycle as a unique attractor in the phase space.  \\
\end{itemize}

To conclude, we observe that $\theta_1 > 0$. If $\theta_1$ was negative, for $\theta=0$ some periodic solution should exist. It is sufficient to show that this is not the case. For $\theta=0$, the trajectory of $\mu$ in \eqref{macrodet} is $\mu(t) = \mu (0) \, e^{-\alpha t}$ and clearly excludes the possibility of having limit cycles.

\subsection{Proof of Proposition \ref{prop:gaussapprox}}

Let $z_{\s}(t)$ be the solution of the linear stochastic differential equation \eqref{linearsde}. Clearly $z_{\s}(\cdot)$ is a centered Gaussian process, so $y_{\s}(t) := m_{\s}(t) + \s z_{\s}(t)$ is obviously Gaussian, which establish property {\bf (a)}. Moreover $\E(y_{\s}(t)) = m_{\s}(t)$, and $Var(y_{\s}(t)) = \s^2 V_{\s}(t)$. Property {\bf (b)}, i.e. the fact that the first two moments of $y_{\s}(t)$ and $\nu_{\s}(\cdot)$ satisfy \eqref{moments} for $k=1,2$, follows from a simple direct computation, that uses equations \eqref{gaussapprox}, and the fact that
\[
\E\left[y^3_{\s}(t) \right] = m^3_{\s}(t) + 3 \s^2 m_{\s}(t)V_{\s}(t), \ \ \ \ \E\left[y^4_{\s}(t) \right] = m^4_{\s}(t) + 6 \s^2 m^2_{\s}(t)V_{\s}(t) + 3 \s^4 V^2_{\s}(t).
\]
So, we are left with the proof of Property {\bf (c)}. By the first equation in \eqref{gaussapprox} and \eqref{linearsde}, we have
\be{yeq}
\begin{split}
dy_{\s} & = \left[ - m_{\s}^3 + m_{\s} - \nu_{\s} - 3 \s^2 m_{\s}V_{\s} - 3 \sigma m_{\s}^2 z_{\s} + \s z_{\s} - 3 \s^3 V_{\s} z_{\s} \right] dt + \s dB \\
& = \left[ - y_{\s}^3 + y_{\s} - \nu_{\s}+ 3 \s^2 m_{\s} (z^2_{\s} - V_{\s}) + \s^3(z_{\s}^3 - 3 V_{\s} z_{\s}) \right] dt + \s dB.
\end{split}
\ee
By subtracting \eqref{yeq} from \eqref{macro}, we get
\be{diffxy}
\begin{split}
x_{\s}(t) - y_{\s}(t) & = \int_0^t \left\{ \left[ y^3_{\s}(s) - x^3_{\s}(s) - y_{\s}(s) + x_{\s}(s) \right] + \left[ \nu_{\s}(s) - \mu_{\s}(s) \right] \right\} ds - \s^2 \int_0^t R_{\s}(s)ds \\
& = \int_0^t \left[x_{\s}(s) - y_{\s}(s)\right] [1-f(s)]ds + \int_0^t \left[ \nu_{\s}(s) - \mu_{\s}(s) \right]ds - \s^2 \int_0^t R_{\s}(s)ds,
\end{split}
\ee
where
\be{remainder}
R_{\s}(t) := 3 m_{\s}(t) (z^2_{\s}(t) - V_{\s}(t)) + \s (z_{\s}^3(t) - 3 V_{\s}(t) z_{\s}(t))
\ee
and
\[
f(t) := x_{\s}^2(t) + y^2_{\s}(t) + x_{\s}(t) y_{\s}(t) \geq 0.
\]
Note that \eqref{diffxy} can be seen as a linear differential equation for $\psi(t) := x_{\s}(t) - y_{\s}(t)$ of the form
\[
\psi(t) = \int_0^t A(s)\psi(s) ds + \int_0^t B(s)ds,
\]
with initial condition $\psi(0) = 0$, whose solution is given by
\be{lineq}
\psi(t) = \int_0^t B(s) \exp\left[\int_s^t A(r) dr \right]ds.
\ee
This gives
\be{x-y1}
x_{\s}(t) - y_{\s}(t) = \int_0^t  e^{-\int_s^t(1-f(r))dr} \left[ \nu_{\s}(s) - \mu_{\s}(s) - \s^2 R_{\s}(s) \right] ds.
\ee
On the other hand, by subtracting from the  equation $\frac{d}{dt} \nu_{\s}(t) = - \a \nu_{\s}(t)  - \theta \frac{d}{dt} \E(y_{\s}(t))$ the corresponding equation for $\mu_{\s}$, we obtain
\[
\nu_{\s}(t) - \mu_{\s}(t) = - \a \int_0^t [\nu_{\s}(s) - \mu_{\s}(s)]ds - \theta[ \E(y_{\s}(t)) -  \E(x_{\s}(t))].
\]
This last formula provides an equation for $\nu_{\s}(t) - \mu_{\s}(t)$ of the form
\[
\psi(t) = \int_0^t A(s)\psi(s) ds + C(t),
\]
whose solution is
\be{lineq1}
\psi(t) = C(t) + \int_0^t A(s) C(s) \exp \left[\int_s^t A(r)dr \right] ds.
\ee
This yields
\[
\nu_{\s}(t) - \mu_{\s}(t) = \theta \left[ \E(x_{\s}(t)) - \E(y_{\s}(t))- \alpha \int_0^t e^{-\a(t-s)} [ \E(x_{\s}(s)) - \E(y_{\s}(s))]ds \right],
\]
from which the following inequality follows:
\begin{align}\label{nu-mu}
\sup_{t \in [0,T]} |\nu_{\s}(t) - \mu_{\s}(t)| &\leq \theta \sup_{t \in [0,T]} \left\vert \E(x_{\s}(t)) - \E(y_{\s}(t)) \right\vert \nonumber\\
& \qquad - \alpha \theta \sup_{t \in [0,T]}  \int_0^t e^{-\a(t-s)} \left\vert \E(x_{\s}(s)) - \E(y_{\s}(s)) \right\vert ds \nonumber\\
&\leq K_T \sup_{t \in [0,T]} \E \left[ \left| x_{\s}(t) - y_{\s}(t) \right| \right],
%
%
%
\end{align}
for some constant $K_T > 0$.
By \eqref{x-y1} and \eqref{nu-mu} we get
\[
\E \left[ \sup_{t \in [0,T]} \left|  x_{\s}(t) - y_{\s}(t) \right| \right] \leq A_T \int_0^T \mathbb{E} \left[ \sup_{s \in [0,t]}  \left|  x_{\s}(s) - y_{\s}(s) \right| \right] ds + A_T \s^2 \sup_{t \in [0,T]}  \E \left[ \left|R_{\s}(t) \right| \right]
\]
for some $A_T>0$ which, by Gronwall's Lemma, yields, for all $\s$ in a bounded subset of $(0,+\infty)$
\be{x-y2}
\E \left[ \sup_{t \in [0,T]} \left|  x_{\s}(t) - y_{\s}(t) \right| \right] \leq B_T \s^2 \sup_{t \in [0,T]}  \E \left[ \left|R_{\s}(t) \right| \right] \leq C_T \s^2,
\ee
where the last inequality follows from the fact that $R_{\s}(t)$ is a polynomial function of a Gauss-Markov process and, therefore, has a $L^1$-norm which is locally bounded in time. To complete the proof of \eqref{normappr} we still need to show that
\[
\sup_{t \in [0,T]} |\mu_{\s}(t) - \nu_{\s}(t) | \leq C_T \s^2
\]
for some $C_T>0$, which follows immediately from  \eqref{x-y2} and \eqref{nu-mu}.

\subsection{Proof of Theorem \ref{th:excitable}}

By symmetry, it is enough to consider the equilibrium $s_1$. Consider the vector field $V_{\a,\theta}^{\s}$ associated to equation \eqref{gaussapprox}
\[
V_{\a,\theta}^{\s}(m,\nu,v) := \left( \begin{array}{c} -m^3 + m - \nu - 3 \s^2 mv \\ -(\a-\theta) \nu + \theta(m^3-m) + 3 \theta \s^2 mv \\ 1+ 2(1-3m^2)v - 6 \s^2 v^2 \end{array} \right).
\]
We linearize the system around $s_1$; we obtain the  linear system associated to the Jacobian matrix
\[
DV_{\alpha,\theta}^{\sigma} \left( s_1 \right) =
\begin{pmatrix}
-1 - \sqrt{1-3\sigma^2} & -1 & 3 \, \sigma^2 \dfrac{\sqrt{1+\sqrt{1-3\sigma^2}}}{\sqrt{2}} \\
\theta \left(1+ \sqrt{1-3\sigma^2}\right) & -\alpha + \theta & -3 \, \theta \sigma^2  \dfrac{\sqrt{1+\sqrt{1-3\sigma^2}}}{\sqrt{2}} \\
\dfrac{3\sqrt{2}}{\sqrt{1+\sqrt{1-3\sigma^2}}}  & 0 & -3 - \sqrt{1-3\sigma^2}
\end{pmatrix}.
\]
In particular, for $\s^2 = 0$ and $\theta = \a+2$,
\[
DV_{\alpha,\alpha+2}^{0} \left( s_1 \right) =
\left(
\begin{array}{cc|c}
-2 & -1 & 0 \\
2 (\alpha + 2) & 2 & 0 \\ \hline
3  & 0 & -4
\end{array}
\right)
\]
whose spectrum is given by 
\[
\mathrm{Spec} \left[ DV_{\alpha,\alpha+2}^{0} \left( s_1 \right) \right] = \{ -4 \} \cup \{ \pm i \sqrt{2 \alpha} \} \,.
\] 
Denote by $\l(\s)$ the continuous function giving an eigenvalue of the matrix $DV_{\alpha,\alpha+2}^{\s} \left( s_1 \right)$, with $\l(0) = \sqrt{2\a} i$. We are going to show that, if $\a<10$, then
\be{posder}
\frac{d}{d(\s^2)} \mathrm{Re}\left(\l(\s) \right)\Big|_{\s=0} >0.
\ee
It follows that there is $\s>0$ such that the matrix $DV_{\alpha,\alpha+2}^{\s} \left( s_1 \right)$ has an eigenvalue with strictly positive real part. By continuity, this implies that if $\theta < \a+2$ but sufficiently close to $\a+2$, then also the matrix $DV_{\alpha, \theta}^{\s} \left( s_1 \right)$ has an eigenvalue with strictly positive real part, and this suffices to complete the proof. Thus, we are left with the proof of \eqref{posder}. We write
\[
0 = \det\left(DV_{\alpha,\alpha+2}^{\s} \left( s_1 \right) - \l(\s) \mathbb{I}  \right) = a_0 (\sigma) + a_1 (\sigma) \lambda (\sigma) + a_2 (\sigma) \left( \lambda (\sigma) \right)^2 + a_3 (\sigma) \left( \lambda (\sigma) \right)^3,
\]
which gives
\[
\frac{d}{d(\s^2)} \left(\l(\s) \right)\Big|_{\s=0} = \frac{1}{a_1(0)} \frac{d}{d(\s^2)} a_0(\s)\Big|_{\s=0}.
\]
Since
\begin{multline*}
 \det\left(DV_{\alpha,\alpha+2}^{\s} \left( s_1 \right) - \l \mathbb{I}  \right) \\ = 
 - \lambda^3 - 2  \left[ 1 + \sqrt{1-3\sigma^2} \right]  \lambda^2 + \left[ 2 - \alpha + 12 \sigma^2 - \left( \alpha + 2 \right) \sqrt{1 - 3\sigma^2} \right] \lambda
\\ - 4 \alpha \sqrt{1 - 3\sigma^2} \left( 1 + \sqrt{1-3\sigma^2} \right),
\end{multline*}
we obtain
\[
\frac{d}{d(\s^2)} \left(\l(\s) \right)\Big|_{\s=0} = \frac{3  \lambda^2(0) + 3 \left( 5 + \frac{\alpha}{2} \right) \lambda(0) + 18 \alpha}{3  \lambda^2(0) + 8 \lambda(0) + 2 \alpha} .
\]
Inserting $\l(0) = \sqrt{2\a} i$, we easily get
\[
\frac{d}{d(\s^2)} \mathrm{Re}\left(\l(\s) \right)\Big|_{\s=0} =\frac{3(10 - \alpha)}{2(8 + \alpha)},
\]
from which \eqref{posder} follows.

\section{Appendix}

\subsection{Proof of Theorem \ref{th:ex}}

It will be convenient to write \eqref{macro} in the form
\be{macroalt}
\begin{split}
dx(t) & =  \left( -x^3(t) + x(t) \right) dt - \mu(t) dt + \s dB(t) \\ 
\frac{d\mu}{dt}(t) & =  - \a \mu(t) - \theta \frac{d}{dt} \E \left[x(t) \right].
\end{split}
\ee
The second equation in \eqref{macroalt} can be solved in $\mu(t)$:
\be{solvemu}
\mu(t) = e^{-\a t}\mu_0 - \theta \E[x(t)] + \theta e^{-\a t} \E(\xi) + \a \theta \int_0^t e^{- \a(t-s)} \E[x(s)]ds.
\ee
To show existence and uniqueness in \eqref{macroalt}  we follow the argument in \cite{Daw83} (Theorem 2.4.1). Given $\mu_0 \in \R$, we define the stochastic processes $(x_n(t))_{t \geq 0}$ and the functions $(\mu_n(t))_{t \geq 0}$, for $n \geq 0$, by the following Picard iteration:
\be{picard}
\begin{split}
dx_{n}(t) & = (-x_n^3(t) + x_n(t))dt - \mu_n(t) dt + \s dB(t) \\
\mu_{n+1}(t) & = e^{-\a t}\mu_0 - \theta \E[x_n(t)] + \theta e^{-\a t} \E(\xi) + \a \theta \int_0^t e^{- \a(t-s)} \E[x_n(s)]ds,
\end{split}
\ee
all with the same initial condition $x_n(0) = \xi$, $\mu_n(0) = \mu_0$.
Setting
\[
f(t) := x^2_{n+1}(t) + x^2_n(t) + x_{n+1}(t) x_n(t) \geq 0,
\]
and observing that
\[
x^3_{n+1}(t) - x^3_n(t) = (x_{n+1}(t) - x_n(t))f(t),
\]
we get
\[
\begin{split}
x_{n+1}(t) - x_n(t) & = \int_0^t \left[ x_{n+1}(s) - x_n(s) \right] [1-f(s)]ds - \int_0^t \left[ \mu_{n+1}(s) - \mu_n(s) \right] ds \\ & = \int_0^t \left[ x_{n+1}(s) - x_n(s) \right] [1-f(s)]ds + \theta \int_0^t \E\left[ x_{n}(s) - x_{n-1}(s) \right]ds \\ & ~~ - \a \theta \int_0^t \int_0^s e^{-\a(s-r)} \E\left[ x_{n}(r) - x_{n-1}(r) \right]drds.
\end{split}
\]
Note that, setting $\psi(t) := x_{n+1}(t) - x_n(t)$, as in the proof of Proposition \ref{prop:gaussapprox}, this last identity is of the form of a linear differential equation
\[
\psi(t) = \int_0^t A(s)\psi(s) ds + \int_0^t B(s)ds,
\]
whose solution is given by \eqref{lineq}.
This yields
\begin{multline}
\label{picard1}
x_{n+1}(t) - x_n(t) \\ = \int_0^t e^{ \int_s^t (1-f(u))du} \left[\theta  \E\left[ x_{n}(s) - x_{n-1}(s) \right] - \a \theta \int_0^s e^{-\a(s-r)} \E\left[ x_{n}(r) - x_{n-1}(r) \right]dr \right] ds.
\end{multline}
Now, set 
\[
\varphi_n(t) := \sup_{s \in [0,t]} \E\left[|x_{n+1}(s) - x_n(s)|\right].
\]
It follows readily from \eqref{picard1} that for each $T>0$ there is a constant $C_T>0$ for which
\[
\varphi_{n+1}(t) \leq C_T \int_0^T \varphi_n(s)ds,
\]
which implies
\[
\varphi_n(T) \leq \varphi_1(T) \frac{(C_T T)^{n-1}}{(n-1)!}.
\]
As a consequence, the sequence of functions $\E[x_n(\cdot)]$ is Cauchy in $\mathcal{C}([0,T])$, and therefore it converges to a limit $m(t)$. Let $x(\cdot)$ be the unique solution of
\be{equnique}
\begin{split}
dx(t) & =  \left( -x^3(t) + x(t) \right) dt - \mu(t) dt + \s dB(t) \\ 
\mu(t) & = e^{-\a t}\mu_0 - \theta m(t) + \theta e^{-\a t} \E(\xi) + \a \theta \int_0^t e^{- \a(t-s)} m(s) ds.
\end{split}
\ee
By mimicking the above argument one gets
\begin{multline}
\label{exist}
x_{n+1}(t) - x(t) \\ = \int_0^t e^{ \int_s^t (1-f(u))du} \left[\theta  \E\left[ x_{n}(s) - x(s) \right] - \a \theta \int_0^s e^{-\a(s-r)} \E\left[ x_{n}(r) - x(r) \right]dr \right] ds,
\end{multline}
where, now, 
\[
f(t) := x^2_{n+1}(t) + x^2(t) + x_{n+1}(t) x(t) \geq 0.
\]
As before, this implies that $\E(x_n(t)) \ra \E(x(t))$, and thus $m(t) = \E(x(t))$. This proves that $x(\cdot)$ is a solution of \eqref{macroalt}. The very same argument is used to prove uniqueness. Indeed, if $y(\cdot)$ is another solution, we write the integral equation for $y(t) - x(t)$ and show as before that $E(x(t)) = E(y(t))$ for all $t$. Thus $y(\cdot)$ and $x(\cdot)$ are both solutions of \eqref{equnique} with the same $m(\cdot)$ and the same initial conditions, from which $x(\cdot) = y(\cdot)$ follows.

\subsection{Proof of Theorem \ref{th:prop-ch}}

We follow a coupling method; we refer to \cite{Szn89} for more details on this approach to prove propagation of chaos. Let $(\xi_i)_{i=1}^N$ be independent random variables with law $\l$, and independent of the Brownian motions $(w_i(t): t \geq 0 )_{i=1}^N$. Moreover, let $(x_i^{(N)}(\cdot), \mu^{(N)}(\cdot))_{i=1}^N$ be the solution of \eqref{micro1} with an initial condition $x^{(N)}_i(0) = \xi_i$, $\mu^{(N)}_i(0) = \mu_0$, while $(y_i(\cdot), \nu(\cdot))_{i=1}^N$ solve \eqref{macro} with the same initial conditions and $B(\cdot) = w_i(\cdot)$.  The claim in Theorem \ref{th:prop-ch} is proved if we show that (the apex $N$ in $x_i^{(N)}$ is omitted in what follows)
\be{convcoupling}
\lim_{N \ra +\infty} \E \left[ \sup_{t \in [0,T]} |x_1(t) - y_1(t)| \right]= 0.
\ee
The strategy is analogous to that of the proof of Theorem \ref{th:ex}. We first write, for $f(t) := x_1^2(t) + y_1^2(t) + x_1(t) y_1(t) \geq 0$:
\[
x_1(t) - y_1(t) = \int_0^t [x_1(s) - y_1(s)][1-f(s)]ds - \int_0^t [\mu(s) - \nu(s)]ds,
\]
yielding
\[
x_1(t) - y_1(t) = \int_0^t [\mu(s) - \nu(s)] \exp\left[ \int_s^t[1-f(r)]dr \right] ds,
\]
which gives, for all $t \in [0,T]$,
\be{chaos1}
|x_1(t) - y_1(t)| \leq C_T \int_0^T | \mu(s) - \nu(s)|ds
\ee
for some constant $C_T>0$. 
In particular, we have
\be{chaos2}
\E\left[ \sup_{t \in [0,T]} |x_1(t) - y_1(t)|  \right] \leq C_T \int_0^T \E[| \mu(s) - \nu(s)|]ds.
\ee

On the other hand, it holds
\begin{multline}\label{chaos3}
\mu(t) - \nu(t) = - \a \int_0^t [\mu(s) - \nu(s)]ds \\
-\theta \left[ \frac{1}{N} \sum_{i=1}^N x_i(t) - \E[y_1(t)] - \left( \frac{1}{N} \sum_{i=1}^N x_i(0) - \E[y_1(0)] \right) \right].
\end{multline}
Set
\[
\varphi(t) := \frac{1}{N} \sum_{i=1}^N x_i(t) - \E[y_1(t)] - \left( \frac{1}{N} \sum_{i=1}^N x_i(0) - \E[y_1(0)] \right).
\]
By \eqref{lineq1},
\be{chaos31}
\mu(t) - \nu(t) = - \theta \varphi(t) + \a \theta \int_0^t \varphi(s) e^{-\a(t-s)}ds.
\ee
Note that if in the expression of $\varphi(t)$ we add and subtract $\frac{1}{N}\sum_i [y_i(t) - y_i(0)]$, the following estimate follows for $s \in [0,T]$:
\begin{align}\label{chaos32}
\E\left[|\varphi(s)|\right] & \leq \frac{1}{N} \sum_{i=1}^N \E[|x_i(s) - y_i(s)|] + \frac{1}{N} \sum_{i=1}^N \E[|x_i(0) - y_i(0)|] \nonumber\\
& \quad + \E \left[ \left| \frac{1}{N} \sum_{i=1}^N y_i(s) - \E[y_1(s)]  \right| \right] + \E \left[ \left| \frac{1}{N} \sum_{i=1}^N y_i(0) - \E[y_1(0)]  \right| \right] \nonumber\\
& \leq \E \left[ \sup_{r \in [0,s]}|x_1(r) - y_1(r)| \right] + \frac{C}{\sqrt{N}},
\end{align}
for some $C>0$,
where the last inequality follows from the fact that  $\E[|x_i(t) - y_i(t)|]$ does not depend on $i$, and that the random variables $y_i(s)$ are i.i.d. and with finite variance.
By \eqref{chaos32} and \eqref{chaos31} it follows that
\be{chaos5}
\E[|\mu(t) - \nu(t)|] \leq A_T \left\{ \E \left[ \sup_{t \in [0,T]} |x_1(t) - y_1(t)| \right] + \frac{C}{\sqrt{N}} \right\}
\ee
for some strictly positive constant $A_T$.
Inserting \eqref{chaos5} in \eqref{chaos2} we obtain, for some possibly different constant $C_T$ and for all $t \in [0,T]$:
\[
\E\left[ \sup_{s \in [0,t]} |x_1(s) - y_1(s)|  \right] \leq C_T \int_0^t \E\left[ \sup_{r \in [0,s]} |x_1(r) - y_1(r)|  \right]ds +  \frac{C_T}{\sqrt{N}},
\]
which implies
\be{chaos6}
\E\left[ \sup_{t \in [0,T]} |x_1(t) - y_1(t)|  \right]  \leq \frac{K_T}{\sqrt{N}}
\ee
for some $K_T>0$, from which the conclusion follows.

\section*{Acknowledgments}

The authors wish to thank Andrea Giacobbe for useful conversations. FC acknowledges financial support of FIRB research grant RBFR10N90W. MF has been partially supported by GA\v{C}R grant P201/12/2613.



\bibliographystyle{plain}

\end{document}